\documentclass[11pt]{amsart}
\usepackage{graphicx,calc,amscd}
\usepackage{epsfig,psfig}
\usepackage{float}
\usepackage{labelfig}
\usepackage{amsmath}
\usepackage{amsfonts}
\usepackage{amssymb}

\newcommand{\CC}{{\mathcal C}}
\newcommand{\BB}{{\mathcal B}}

\newcommand{\Sp}{{\mathbb S}}
\newcommand{\G}{{\mathcal G}}

\newcommand{\ii}{\mathrm{int}}

\renewcommand{\tt}{\mbox{\rm Tr}_{2^d|2}\,}
\newcommand{\arcsinh}{\mathrm{arcsinh}}
\newcommand{\arccosh}{\mathrm{arccosh}}

\newtheorem{theorem}{Theorem}[section]
\newtheorem{lemma}[theorem]{Lemma}

\newtheorem{proposition}[theorem]{Proposition}
\newtheorem{remark}[theorem]{Remark}

\begingroup\makeatletter\ifx\SetFigFont\undefined%
\gdef\SetFigFont#1#2#3#4#5{%
 \reset@font\fontsize{#1}{#2pt}%
 \fontfamily{#3}\fontseries{#4}\fontshape{#5}%
 \selectfont}%
\fi\makeatother\endgroup%
\begin{document}
\title{Fixed point free involutions on Riemann surfaces}
\author{Hugo Parlier}
\address{Departamento de Matem\'{a}ticas Fundamentales\\
Facultad de Ciencias\\
Universidad Nacional de Educaci\'{o}n a Distancia\\
Madrid 28040\\
SPAIN} \email{hugo.parlier@epfl.ch}

\thanks{The author was supported in part by the Swiss National Science Foundation
grants 20 - 68181.02 and PBEL2-106180.} 
\date{\today}
\begin{abstract}Involutions without fixed points on hyperbolic closed Riemann
surface are discussed. For an orientable surface $X$ of even genus
with an arbitrary Riemannian metric $d$ admitting an involution
$\tau$, it is known that $\min_{p\in X}d(p,\tau(p))$ is bounded by
a constant which depends on the genus of $X$. The equivalent
result is proved to be false in odd genus, and the optimal
constant for hyperbolic Riemann surfaces is calculated in genus
$2$.
\end{abstract}
\maketitle

\section{Introduction}                  \label{Sect:S1}

Involutions play an important role in the study of compact Riemann
surfaces. For instance, the study of hyperelliptic surfaces of
genus $g$ corresponds to the study of surfaces admitting an
orientation preserving involution with $2g+2$ fixed points and a
Klein surface is the quotient of a Riemann surface by an
orientation reversing involution. Furthermore, by the
uniformization theorem, surfaces with an orientation reversing
involution are conformally equivalent to real algebraic curves. A
further motivation can be found in \cite{bacrivka05} where the
so-called area filling conjecture is treated. The conjecture,
first found in \cite{gr83} for the $n$-dimensional case, states
(in two dimensions) that any surface $S$ with a simple boundary,
endowed with a Riemannian metric with the property that
diametrally opposite points on the boundary are of distance
$\geq\pi$, has area greater or equal to $2\pi$, equality occurring
only in the case of the classical
hemisphere. The conjecture is equivalent to the following:\\

\noindent{\bf Conjecture}\hspace{0.3cm} Let $S$ be an orientable
surface of even genus with a Riemannian metric $d$ that admits a
fixed point free, orientation reversing involution $\tau$. Then
there is a point $p\in S$ with

\begin{equation}\label{eqn:conj}
\frac{d(p,\tau(p))^2}{{\rm area}({\mathcal S})} \leq
\frac{\pi}{4}.
\end{equation}

In \cite{bacrivka05}, this problem is solved when $S$ is
hyperelliptic, thus in particular for genus $2$. Furthermore, if
$S$ is not hyperelliptic, in \cite{ivka04} it is shown that

\begin{equation}\label{eqn:hypevengenus} \frac{d(p,\tau(p))^2}{{\rm area}({\mathcal S})} \leq C
\end{equation}

where $C\in [\frac{\pi}{4},1]$. If $S$ is a hyperbolic Riemann
surface of even genus $g$, this shows that

\begin{equation}\label{eqn:hypRS} d(p,\tau(p))\leq \sqrt{C 2\pi(2g-2)}.
\end{equation}

The idea of the paper is to treat these questions in more detail
for hyperbolic Riemann surfaces. The first main result of the
article is to show that the conjecture stated above cannot be
extended to odd genus, even if one restricts oneself to
hyperelliptic surfaces.

\begin{theorem} For any odd $g\geq 3$ and positive constant $k$, there exists a
hyperbolic Riemann surface $S$ of genus $g$ admitting an
orientation reversing involution $\tau$ which verifies
$d(p,\tau(p))>k$ for all $p\in S$. Furthermore, $S$ can be chosen
hyperelliptic.
\end{theorem}

The corresponding problem for orientation preserving involutions
is also treated, with the same result.\\

The bound in equation \ref{eqn:hypRS} is not sharp and although
the bound in \ref{eqn:hypevengenus} is sharp for arbitrary
metrics, it is not sharp for hyperbolic metrics. The second main
result of the article is the sharp bound in genus $2$ for
hyperbolic metrics.

\begin{theorem} Let $S$ be a Riemann surface of genus $2$ endowed with a
hyperbolic metric and a fixed point free involution $\tau$. On
$S$, there is a $p$ such that $d(p,\tau(p))\leq
\arccosh(\frac{5+\sqrt{17}}{2})$. This upper-bound is sharp and is
attained by a unique hyperbolic Riemann surface of genus $2$ (up
to isometry).
\end{theorem}

Compare the optimal bound ($=2.19...$ ) with the bound from
equation \ref{eqn:hypRS} ($=3.54...$) applied to hyperbolic
surfaces. The unique surface which attains this upper bound is
surprisingly not the surface known as the Bolza curve, which is
maximal for systole length, number of systoles ($12$), number of
automorphisms ($48$), and verifies an optimal systolic inequality
for hyperelliptic invariant CAT(0) metrics in genus $2$
(see \cite{sc931} and \cite{kasa05}).\\

\section{Definitions, notations and preliminaries}
In this article, a {\it surface} will always mean a compact
Riemann surface endowed with a hyperbolic metric. We shall suppose
that the boundary of a surface is a collection of simple closed
geodesics. The signature of the surface will be denoted $(g,n)$
(genus and number of boundary components) and verifies $(g,n) \geq
(0,3)$ (with lexicographic ordering) and $(g,n)\neq (1,0)$. This
condition is imposed by the existence of a hyperbolic metric. A
surface of signature $(0,3)$ is commonly called a $Y$-piece, or a
pair of pants. The area of the surface is given by ${\rm area}(S)
= -2\pi \chi(S)$ where $\chi(S)=2-2g -n$ is the Euler
characteristic of the surface. Unless mentioned, a {\it geodesic}
is a simple closed geodesic on $S$. Distance on $S$ (between
subsets, points or curves) is denoted $d(\cdot,\cdot)$. Curves and
geodesics will be considered primitive and non-oriented, and can
thus be seen as point sets on $S$. A geodesic $\gamma$ (resp. a
set of geodesics $E$) is called {\it separating} if $S\setminus
\gamma$ (resp. $S\setminus E$) is not connected. Let us recall
that a non-trivial closed curve $c$ (not necessarily simple) on
$S$ is freely homotopic to exactly one closed geodesic $\gamma$.
(We denote $\G(c)=\gamma$). If $c$ is simple, then so is $\gamma$.
We shall denote the intersection number of two geodesics $\gamma$
and $\delta$ by $\ii(\gamma,\delta)$. The length of a path or a
curve will be denoted $\ell(\cdot)$, although in general a curves
name and its length will not be distinguished. The {\it systole}
$\sigma$ of a surface is the (or a) shortest non-trivial closed
curve on $S$ (sometimes called a systolic loop, although for
Riemann surfaces, systole seems to be the standard denomination).
A systole is always a simple closed geodesic and cannot intersect
another systole more than once. An involution $\tau$ is an
isometric automorphism of the surface that is of order $2$.
Involutions can either be orientation preserving or not. For
Riemann surfaces this is equivalent to whether the involution is
holomorphic or antiholomorphic. Let us also recall that two
geodesics of length less or equal to $2\arcsinh 1$ are disjoint
(and simple). This can be seen using the following well known
result, commonly called the collar theorem (i.e.
\cite{kee74}, \cite{bu78}, \cite{bubook}, \cite{ra79}).\\

\begin{theorem} Let $\gamma_1$ and $\gamma_2$ be non-intersecting
simple closed geodesics on $S$. Then the collars
$$
\CC(\gamma_i)=\{p\in S \mid d_S(p,\gamma_i)\leq w(\gamma_i)\}
$$
of widths
$$
w(\gamma_i)=\arcsinh(1/\sinh\frac{\gamma_i}{2})
$$
are pairwise disjoint for $i=1,2$. Furthermore, each
$\CC(\gamma_i)$ is isometric to the cylinder
$[-w(\gamma_i),w(\gamma_i)]\times \Sp^1$ with the metric
$ds^2=d\rho^2+\gamma_i^2\cosh^2\!\rho \, dt^2$.
\end{theorem}

Notice that the $\gamma_i$'s divide their collar into two
connected spaces which we will call {\it half-collars}. In the
sequel we will make use of the fact that the collars of two
disjoint geodesics are also disjoint.\\

Simple closed geodesics and orientation reversing involutions are
closely related. The following proposition is an extension of what
is generally called Harnack's theorem \cite{wephd} and can be
found in \cite{krma82}.

\begin{proposition}\label{prop:harnackgen}
If a surface $S$ admits $\tau$, an orientation reversing
involution, then the fixed point set of $\tau$ is a set of $n$
disjoint simple closed geodesics
${\BB}=\{\beta_1,\hdots,\beta_n\}$ with $n\leq g+1$. In the case
where the set $\BB$ is separating, then $S\setminus \BB$ consists
of two connected components $S_1$ and $S_2$ such that $\partial
S_1 =\partial S_2 = \BB$ and $S_2=\sigma(S_1)$. If not, then $\BB$
can be completed by either a set $\alpha$ which consists of one or
two simple closed geodesics such that $\BB \cup \alpha$ has the
properties described above (with the important difference that
$\alpha$ does not contain any fixed points of $\sigma$). Each of
the simple closed geodesics in $\alpha$ is globally fixed by
$\sigma$.
\end{proposition}

\section{The general case}
Let $S$ be a surface of genus $g\geq 2$, endowed with a hyperbolic
metric and a fixed point free involution
$\tau$.\\

\begin{proposition}\label{prop:evenOR} If $g$ is even, then $\tau$ is orientation
reversing. \end{proposition}

\begin{proof} The proof is a consequence of the Riemann-Hurwitz formula,
or can be seen as follows. In general, if $\tau$ is an orientation
preserving isometry, then $S / \tau= O$ is an orientable orbifold
with singular points who lift on $S$ to fixed points of $\tau$. As
$\tau$ is without fixed points, $O$ is an orientable closed
surface of genus $g_o$ endowed with a hyperbolic metric and thus
of area $\frac{{\rm area}(S)}{2}$. This implies that $O$'s Euler
characteristic is equal to $(g-1)$. From this we have
$g_o=\frac{g-1}{2}$ which is not possible if $g$ is even.
\end{proof}

\begin{remark}
As the essence of the proof is purely topological, the proposition
holds for arbitrary metric.
\end{remark}

The following propositions concern further relationships between
simple closed geodesics and involutions. Suppose that $S$ is a
hyperbolic surface of genus $g\geq 2$ and $\tau$ is a fixed point
free involution.

\begin{proposition} Let $\sigma$ be a systole of $S$. Then $\sigma
\cap \tau(\sigma) = \emptyset$ or $\sigma=\tau(\sigma)$.
\end{proposition}

\begin{proof} Suppose that $\sigma \neq \tau(\sigma)$. As $\tau(\sigma)$ is necessarily another systole of
$S$, then the curves $\sigma$ and $\tau(\sigma)$ cannot intersect
more than once (this would imply the existence of a shorter
non-trivial closed curve on $S$). Suppose that $\sigma\cap
\tau(\sigma)= p$, where $p$ is a point. Then $\tau(p)=p$ which
contradicts the hypotheses.
\end{proof}

As $S$ is compact, the value $\min_{p\in S}d(p,\tau(p))$
exists and is attained for at least one point.\\

\begin{proposition}\label{prop:minonsimple}
Let $p\in S$ such that $d(p,\tau(p))$ is minimum. Then $p$ lies on
a simple closed geodesic $\gamma$ of length $2d(p,\tau(p))$ that
verifies $\gamma=\tau(\gamma)$.
\end{proposition}

\begin{proof}
Let $p$ be such a point. Let $c_p$ be a minimal path between $p$
and $\tau(p)$ (thus $d(p,\tau(p)) = \ell(c_p)$). Notice that
$c_p\cap \tau(c_p)=\{p,\tau(p)\}$, and the two paths $c_p$ and
$\tau(c_p)$ are not freely homotopic among simple paths with
endpoints on $p$ and $\tau(p)$. Thus $c_p\cup \tau(c_p)$ is a
non-trivial simple closed curve. Furthermore, it follows that
$\gamma=\tau(\G(c_p\cup \tau(c_p))=\G(c_p\cup \tau(c_p))$ and thus
for $q\in\gamma$, $d(q,\tau(q))\leq d(p,\tau(p))$ with equality
occurring only in the case where $\G(c_p\cup \tau(c_p))=c_p\cup
\tau(c_p)$. This implies that $p$ lies on $\gamma$ and that
$\gamma$ is of length $2d(p,\tau(p))$.
\end{proof}

Notice that for such a geodesic $\gamma$, all point are
diametrically opposite to their images by $\tau$. In fact this is
true for any simple geodesic left invariant by $\tau$.

\begin{proposition} Let $\gamma$ be a simple closed geodesic such that
$\tau(\gamma)=\gamma$. Then the image $\tau(p)$ of $p\in \gamma$
is the point on $\gamma$ diametrically opposite from $p$.
\end{proposition}

\begin{proof} For any $p\in \gamma$, the points $p$ and $\tau(p)$
separate $\gamma$ into two geodesic arcs. The image of one of the
arcs has to be the other arc (otherwise $\tau$ has fixed points)
and since $\tau$ is an isometry, the result follows.
\end{proof}

As mentioned in the introduction, a result in \cite{ivka04}
implies that for any $S$ of even genus with an involution $\tau$
(necessarily orientation reversing by proposition
\ref{prop:evenOR})

\begin{equation*} d(p,\tau(p))\leq \sqrt{C 2\pi(2g-2)}
\end{equation*}

with $C\in[\frac{\pi}{2},2]$. This is false in odd genus for both
orientation preserving and reversing involutions.

\begin{theorem} Let $g\geq 3$ be an odd integer and $k$ a positive constant. There exists a
hyperbolic Riemann surface $S$ of genus $g$ admitting an
orientation preserving involution $\tau$ which verifies
$d(p,\tau(p))>k$ for all $p\in S$. The same result holds for
orientation reversing involutions. Furthermore, in both cases, $S$
can be chosen hyperelliptic.
\end{theorem}

\begin{proof}
We shall begin by showing the general idea for constructing
surfaces with orientation preserving involutions, then surfaces
with orientation reversing involutions which verify the conditions
of the theorem. Finally, we shall give an explicit example of a
hyperelliptic surface which has both an orientation preserving
involution, and an orientation reversing involution which verifies
the conditions of the theorem.\\

Consider a surface $\tilde{S}$ of signature $(\tilde{g},2)$ with
boundary geodesics $\alpha$ and $\beta$ of equal length $x$. Let
$p_\alpha$ and $q_\alpha$ be two diametrically opposite points on
$\alpha$, and $p_\beta$ and $q_\beta$ be two diametrically
opposite points on $\beta$. Consider two copies of $\tilde{S}$,
say $\tilde{S}_1$ and $\tilde{S}_2$, and paste them along their
boundary geodesics such that $\tilde{S}_1$'s $\alpha$ is pasted to
$\tilde{S}_2$'s $\beta$, and $\tilde{S}_2$'s $\alpha$ is pasted to
$\tilde{S}_1$'s $\beta$. The pasting must respect the choice of
$p$s and $q$s, meaning $p_\alpha$ is pasted to $p_\beta$ etc. The
resulting surface $S$ is an orientable surface of genus $2g+1$
admitting an orientation preserving involution $\tau_o$ without
fixed points which acts as follows: a point originally on
$\tilde{S}_1$ is sent to it's corresponding point on $\tilde{S}_2$
and vice-versa. The simple closed geodesics of $S$, previously the
boundary geodesics of $\tilde{S}_1$ and $\tilde{S}_2$, are of
length $x$ and are reversed by $\tau_o$.

\begin{figure}[h]
\leavevmode \SetLabels
\L(.46*.89) $\alpha_1$\\
\L(.575*.08) $\alpha_2$\\
\L(.46*.08) $\beta_1$\\
\L(.575*.89) $\beta_2$\\
\L(.33*.4) $\tilde{S}_1$\\
\L(.65*.4) $\tilde{S}_2$\\
\endSetLabels
\begin{center}
\AffixLabels{\centerline{\epsfig{file = 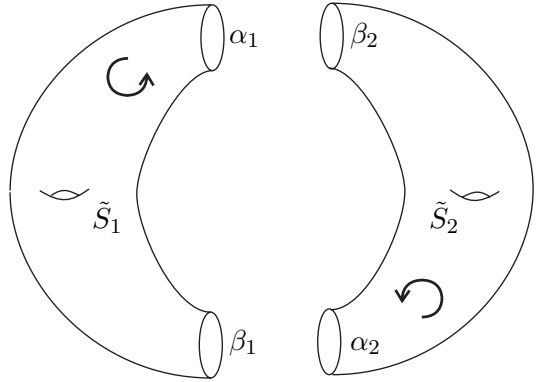,
width=7cm,angle= 0}}}
\end{center}
\caption{Example in genus $3$ with orientation preserving
involution} \label{fig:OP}
\end{figure}

The image $\tau_o(p)$ of a point $p\in S$ is at least ``half a
collar away" from $p$, and by the collar theorem the following
inequality is thus verified:
$$
d(p,\tau_o(p))> \arcsinh(1/\cosh(\frac{x}{2})).
$$

The half-collar length tends to infinity as $x$ tends to $0$, thus
for any $k>0$, it suffices to chose $x$ such that $k <
\arcsinh(1/\cosh(\frac{x}{2})$, and the result follows.\\


Now let us treat the case of orientation reversing involutions.
Consider a $\tilde{S}$ as above. Instead of pasting two identical
copies of $\tilde{S}$, consider $\tilde{S}$ and a symmetric copy
of $\tilde{S}$ (a mirror image), say $\tilde{S}_-$. Denote by
$\alpha_-$ and $\beta_-$ the images of $\alpha$ and $\beta$ on
$\tilde{S}_-$ as in the following figure.

\begin{figure}[h]
\leavevmode \SetLabels
\L(.46*.89) $\alpha$\\
\L(.575*.08) $\beta_-$\\
\L(.46*.08) $\beta$\\
\L(.575*.89) $\alpha_-$\\
\L(.33*.4) $\tilde{S}$\\
\L(.65*.4) $\tilde{S}_-$\\
\endSetLabels
\begin{center}
\AffixLabels{\centerline{\epsfig{file = 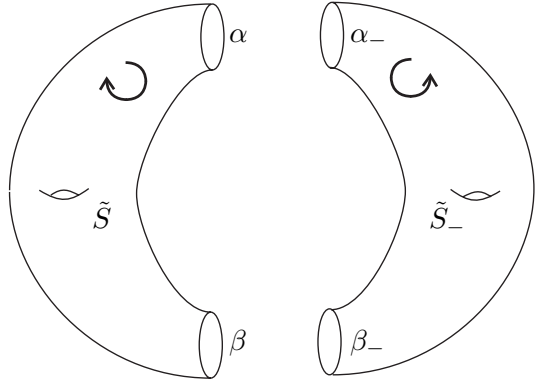,
width=7cm,angle= 0}}}
\end{center}
\caption{Example in genus $3$ with orientation reversing
involution} \label{fig:OR}
\end{figure}

Paste the boundary geodesics together ($\alpha$ to $\alpha_-$,
$\beta$ to $\beta_-$) while respecting the choice of $p$s and
$q$s. The resulting surface $S$ is of genus $2g+1$ and the
orientation reversing involution $\tau_r$ is the isometry taking a
point of $\tilde{S}$ to its corresponding point on $\tilde{S}_-$
and vice-versa. As in the first examples, for any $k$, $x$, the
length of both $\alpha$ and $\beta$, can be chosen such that
$d(p,\tau_r(p))>k$ for all $p\in S$.\\

We shall now give an explicit example of a hyperelliptic surface
with both an orientation preserving involution and an orientation
reversing involution which verifies the conditions of the
theorem.\\

Consider a right-angled $2g+4$-gon in the hyperbolic plane, say
$P$, with edges labeled in cyclic ordering
$\{a_1,b_1,\hdots,a_{g+2},b_{g+2}\}$. The figures are all done
when $g=3$.

\begin{figure}[H]
\leavevmode \SetLabels
\L(1.01*.91) $a_1$\\
\L(.5*1.03) $b_1$\\
\L(.03*-.08) $a_2$\\
\L(.13*.27) $b_2$\\
\L(.23*-.08) $a_3$\\
\L(.33*.27) $b_3$\\
\L(.43*-.08) $a_4$\\
\L(.53*.27) $b_4$\\
\L(.615*-.08) $a_5$\\
\L(.75*.57) $b_5$\\
\endSetLabels
\begin{center}
\AffixLabels{\centerline{\epsfig{file = 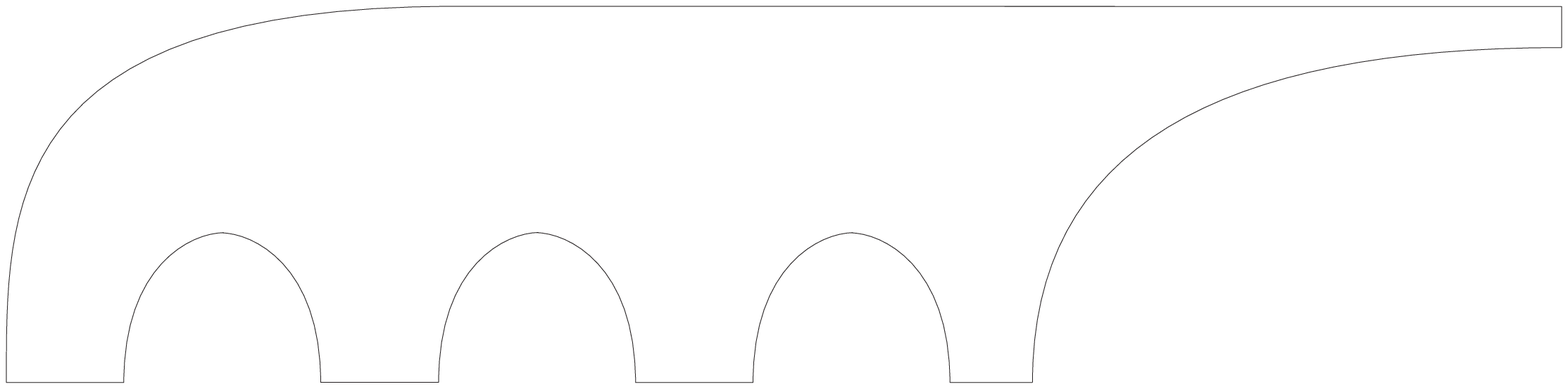,
width=14cm,angle= 0}}}
\end{center}
\caption{The polygon $P$} \label{fig:tengon}
\end{figure}

Take two isometric copies of $P$, say $P_1$ and $Q_1$, and two
symmetric images of $P$, say $P_2$ and $Q_2$ and paste them along
the sides labeled $a_i$ as in the following figure to obtain a
surface $S^+$ of signature $(0,2g+2)$. The bold curves represent
the boundary geodesics of $S^+$.

\begin{figure}[h]
\leavevmode \SetLabels
\L(.2*.8) $P_1$\\
\L(.2*.15) $P_2$\\
\L(.76*.15) $Q_1$\\
\L(.76*.8) $Q_2$\\
\endSetLabels
\begin{center}
\AffixLabels{\centerline{\epsfig{file = 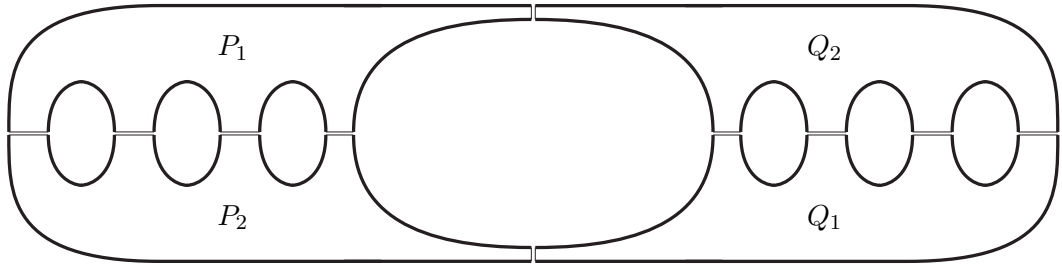,
width=14cm,angle= 0}}}
\end{center}
\caption{The surface $S^+$} \label{fig:polygonsurface}
\end{figure}

By taking a symmetric copy $S^-$ of the surface thus obtained, and
denote by $P_1^-$ etc. the various images of the polygons of
$S^+$. By pasting the two surfaces $S^+$ and $S^-$ along the
$2g+2$ boundary geodesics, as in the following figure, one obtains
a surface $S$ of genus $2g+1$. We require the pasting to be exact
- each edge of a polygon is pasted exactly to its corresponding
edge with end points of each edge coinciding.\\

\begin{figure}[h]
\leavevmode \SetLabels
\L(.2*.63) $P_1$\\
\L(.2*.12) $P_2$\\
\L(.76*.12) $Q_1$\\
\L(.76*.63) $Q_2$\\
\L(.23*1.035) $P_1^-$\\
\L(.41*.26) $P_2^-$\\
\L(.52*.23) $Q_1^-$\\
\L(.675*1.035) $Q_2^-$\\
\endSetLabels
\begin{center}
\AffixLabels{\centerline{\epsfig{file = 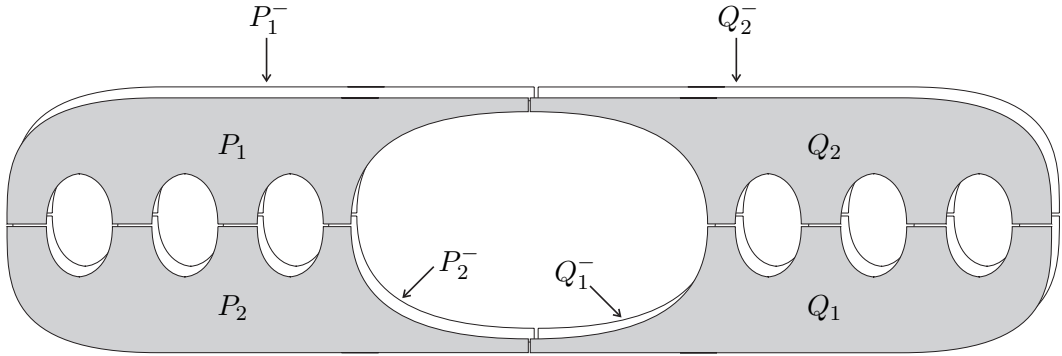,
width=14cm,angle= 0}}}
\end{center}
\caption{The surface $S$} \label{fig:polygonsurfacetotal}
\end{figure}

The surface thus obtained necessarily admits a number of
involutions. One of these is the hyperelliptic involution
$\tau_h$, with fixed points the end points of all the images of
$a_3,\hdots,a_{g+1}$. (To be precise, $\tau_h$ exchanges $P_1$ and
$P_2^-$, $P_2$ and $P_1^-$, $Q_1$ and $Q_2^-$, $Q_2$ and
$Q_1^-$.)\\

The remaining involutions we are interested in are the ones
without fixed points. The orientation preserving involution
$\tau_o$ defined by exchanging $P_1$ with $Q_1$, $P_2$ with $Q_2$,
$P_1^-$ with $Q_1^-$, and $P_2^-$ with $Q_1^-$, does not have any
fixed points. The orientation reversing involution $\tau_r$
defined by exchanging $P_1$ with $Q_1^-$, $P_2$ with $Q_2^-$,
$P_1^-$ with $Q_1$, and $P_2^-$ with $Q_1$, does not have any
fixed points either. The eight images of the edge $a_1$ on $S$
form two simple closed geodesics (after pasting) of length $2a_1$.
As in the more general case, described above, for any $k$, we can
chose a sufficiently short $a_1$ such that we have
$d(p,\tau_o(p))>k$ as well as $d(p,\tau_r(p))>k$ for all $p\in S$.
\end{proof}

\section{The case of genus $2$}

By the previous theorem, for all surfaces of genus $2$, the
minimum distance between a point and the image by an involution is
bounded. The object of this section is a detailed study of fixed
point free involutions on surfaces of genus $2$ leading to the
sharp bound on this distance.

\begin{lemma}\label{lemma:genus21} Let $\tau$ be an involution without fixed points on a surface of
genus $2$. Then the following statements are true:\\
\begin{enumerate}
\item $\tau$ reverses orientation.\\
\item $S$ contains a separating simple closed geodesic $\beta$
such that $\beta=\tau(\beta)$. The two parts of $S$ separated by
$\beta$ are interchanged by $\tau$.\\
\item For any $p\in \beta$, $p$ and $\tau(p)$ are diametrically
opposite.\\
\end{enumerate}
\end{lemma}

\begin{proof}
Parts (1) and (3) have been treated earlier in the general case
and part (2) is a direct consequence of proposition
\ref{prop:harnackgen} (Harnack's theorem). \end{proof}

The geodesic $\beta$ from the previous lemma divides $S$ into two
surfaces of signature $(1,1)$. The following lemma recalls some
essential facts about these surfaces. These facts are either well
known, or their proofs can be found in \cite{sc931} (theorem 4.2,
p. 578). For such surfaces, and a choice of interior simple closed
geodesic $\alpha$, we denote $h_\alpha$ the unique simple geodesic
path which goes from boundary to boundary and intersects boundary
at two right angles and does not cross $\alpha$. We will refer to
the geodesic path $h_\alpha$ as the height associated to $\alpha$.

\begin{figure}[h]
\leavevmode \SetLabels
\L(.466*.782) $\alpha$\\
\L(.457*.38) $h_\alpha$\\
\endSetLabels
\begin{center}
\AffixLabels{\centerline{\epsfig{file =
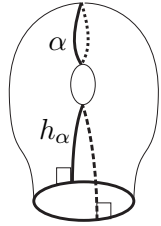,width=2cm,angle=0}}}
\end{center}
\caption{The curve $h_\alpha$ associated to $\alpha$}
\label{fig:qpieceheight}
\end{figure}

\begin{lemma}\label{lemma:genus22} Let $Q$ be a surface of signature $(1,1)$ with boundary geodesic $\beta$. Then the
following statements are true:

\begin{enumerate}
\item $Q$ is hyperelliptic and its hyperelliptic involution has
three fixed points in the interior of $Q$ (the Weierstrass points
of $Q$).\\

\item Let $\gamma$ be an interior simple closed geodesic of $Q$ and denote its
associated height $h_\gamma$. Then $\gamma$ passes through exactly
two of the three Weierstrass points and the remaining Weierstrass
point is the midpoint of $h_\gamma$. Furthermore, the length of
$\gamma$ is directly proportional to the length of $h_\gamma$.\\

\item
Among all surfaces of boundary length $\ell(\beta)$, the unique
surface (up to isometry) with maximum length systole is the
surface with three distinct systoles.

\item
If $Q$ has three systoles, then their associated heights do not
intersect and are evenly spaced along $\beta$. \end{enumerate}
\end{lemma}

We shall now proceed to the main result of this section.

\begin{theorem} Let $S$ be a Riemann surface of genus $2$ endowed with a
hyperbolic metric and a fixed point free involution $\tau$. On
$S$, there is a $p$ such that $d(p,\tau(p))\leq
\arccosh(\frac{5+\sqrt{17}}{2})$. This upper-bound is sharp and is
attained by a unique hyperbolic Riemann surface of genus $2$ (up
to isometry).
\end{theorem}

\begin{proof}
For $S$ of genus $2$, let $\beta$ be the separating geodesic
described in the previous lemma \ref{lemma:genus21}. Let $Q$ be
one of the two surfaces of signature $(1,1)$ separated by $\beta$.
Proposition \ref{prop:minonsimple} implies that a $p$ which
minimizes $d(p,\tau(p)$ is found on a simple closed geodesic left
invariant by $\tau$. As all simple closed geodesics in the
interior of both $Q$ and $\tau(Q)$ are distinct from their images,
$p$ must be found on a simple geodesic that either crosses $\beta$
or is $\beta$.\\

Let $h$ be a height on $Q$. Then $h\cup \tau(h)=\gamma_h$ is a
simple closed geodesic with the property that for $p\in \gamma_h$,
$\tau(p)$ is diametrically opposite to $p$ (lemma
\ref{lemma:genus21}). Furthermore, if $\sigma$ is a systole on
$Q$, then $h_\sigma$ is the shortest height on $Q$. Notice that
all simple closed geodesics that cross $\beta$ are longer that
$2h_\sigma$, thus $d(p,\tau(p)) \geq h_\sigma$, equality occurring
when $2h_\sigma\leq \beta$. For any $\beta$, the maximum value
that $h_\sigma$ can attain is attained in the situation described
in lemma \ref{lemma:genus22}. Furthermore, the maximum value of
$h_\sigma$ is strictly proportional to the length of $\beta$. The
maximum value for $d(p,\tau(p))$ is thus obtained when both
$\beta/2$ and maximum $h_\sigma$ are equal. These conditions
define a unique surface of signature $(1,1)$, whose lengths verify
(i.e. \cite{sc931})

\begin{equation}\label{eqn:proofthmB1}
4\cosh^3(\frac{\sigma}{2})+6\cosh^2(\frac{\sigma}{2})+1=\cosh(\frac{\beta}{2}).
\end{equation}

\begin{figure}[h]
\leavevmode \SetLabels
\L(.305*.75) $\frac{\sigma}{2}$\\
\L(.45*.41) $\frac{h_\sigma}{2}$\\
\L(.38*.46) $\mathbf P$\\
\L(.40*.09) $\frac{\beta}{4}$\\
\endSetLabels
\begin{center}
\AffixLabels{\centerline{\epsfig{file =
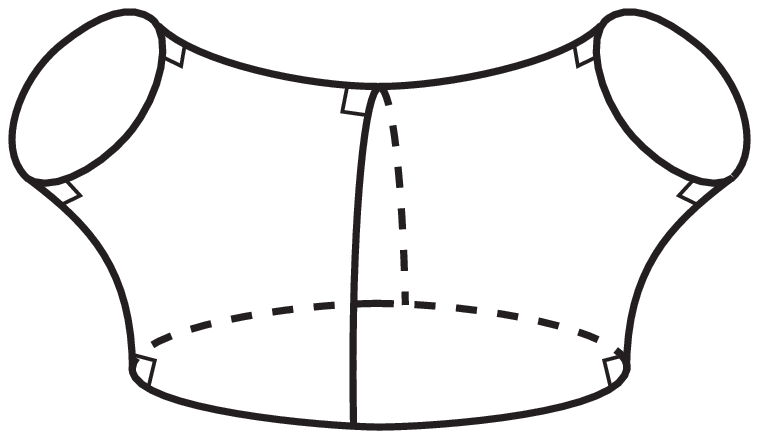,width=7cm,angle=0}}}
\end{center}
\caption{Maximal $Q$ cut along a systole} \label{fig:proofthmB}
\end{figure}

The figure shows this surface after cutting along a systole
$\sigma$. The hyperbolic right angled pentagon $\mathbf P$ has
equal length adjacent edges of lengths $\frac{h_\sigma}{2}$, and
$\frac{\beta}{4}$, and opposite edge of length $\frac{\sigma}{2}$.
Using the hyperbolic trigonometry formula for such pentagons, one
obtains

$$
\sinh^2(\frac{\beta}{4})= \cosh(\frac{\sigma}{2})
$$

and thus

\begin{equation}\label{eqn:proofthmB2}
\cosh(\frac{\beta}{2})-1= 2 \cosh(\frac{\sigma}{2}).
\end{equation}

Using equations \ref{eqn:proofthmB1} and \ref{eqn:proofthmB2}, the
value for a systole $\sigma$ verifies

\begin{equation}\label{eqn:proofthmB3}
2\cosh^2(\frac{\sigma}{2})-3\cosh(\frac{\sigma}{2})-1=0.
\end{equation}

From this we can deduce the value

$$h_\sigma=\frac{\beta}{2}=\arccosh(\frac{5+\sqrt{17}}{2}).$$

In order to obtain the maximal surface of genus $2$, it suffices
to take two copies of the maximum surface of signature $(1,1)$ and
to paste them such that minimum length heights touch. As the
heights are evenly spaced along $\beta$, whichever way this is
done one will obtain the same surface. This surface, say
$S_{\max}$, is clearly unique up to isometry. The following figure
explicitly illustrates how to obtain $S_{\max}$ by pasting $8$
copies of $\mathbf P$. On $\mathbf P$, consider the midpoint of
the edge labeled $\frac{\sigma}{2}$ on figure \ref{fig:proofthmB}.
The points labeled $p_1$, $p_2$, $q_1$ and $q_2$ are the $8$
copies of this point and the points $p_3$ and $q_3$ are as labeled
on the figure. The pasting of the boundary geodesics is as
indicated in the figure ($p_1$ pasted to $p_1$ etc.). The
Weierstrass points of $S_{\max}$ are exactly the points
$p_1,p_2,p_3$ and $q_1,q_2,q_3$.

\begin{figure}[h]
\leavevmode \SetLabels
\L(.33*.77) $p_1$\\
\L(.25*.91) $p_2$\\
\L(.65*.77) $p_2$\\
\L(.73*.91) $p_1$\\
\L(.33*.21) $q_1$\\
\L(.25*.07) $q_2$\\
\L(.65*.21) $q_2$\\
\L(.73*.07) $q_1$\\
\L(.49*.94) $p_3$\\
\L(.49*.05) $q_3$\\
\endSetLabels
\begin{center}
\AffixLabels{\centerline{\epsfig{file =
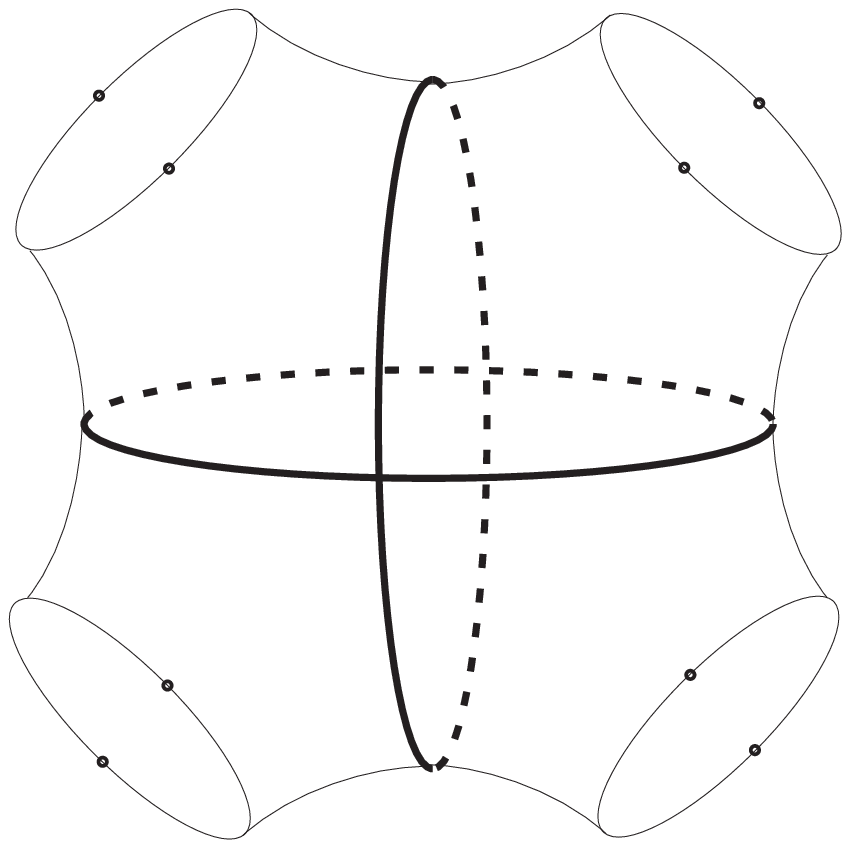,width=8cm,angle=0}}}
\end{center}
\caption{$S_{\max}$} \label{fig:smax}
\end{figure}

Notice that their are exactly $6$ systoles on this surface, and
the length of a systole of the surface is the same as a systole on
the maximal $(1,1)$ surface. The length of a systole, by equation
\ref{eqn:proofthmB3} is

$$
\sigma=2\arccosh(\frac{3+\sqrt{17}}{4}).
$$
\end{proof}

The maximal surface $S_{\max}$ we have constructed above has a
remarkable property: it is {\it not} the Bolza curve. The Bolza
curve, with its unique hyperbolic metric, is constructed in the
same fashion, namely it is obtained by pasting two maximal
surfaces of signature $(1,1)$. It is the unique maximal genus $2$
surface for systole length. Compared to $S_{\max}$, its separating
geodesic $\beta$ is longer and its systole length is
$2\arccosh(1+\sqrt{2})$.

\bibliographystyle{plain}
\def\cprime{$'$}

\end{document}